\documentclass[letterpaper,journal,twocolumn,twoside,final,9pt]{IEEEtran}
 \usepackage{amssymb}
 \usepackage{graphicx}

\providecommand\href[2]{#2}

\newtheorem{theorem}{Theorem}[section]
\newtheorem{proposition}[theorem]{Proposition}
\newtheorem{corollary}[theorem]{Corollary}
\newtheorem{lemma}[theorem]{Lemma}
\newtheorem{example}[theorem]{Example}

\renewcommand{\title}[1]{\begin{center}{\large\textbf{#1}}\\[2.0ex]}
\renewcommand{\author}[1]{{\rm #1}\end{center}\vspace{2ex}}
\renewcommand{\maketitle}{}
\renewcommand{\thanks}{\footnotetext}

\begin{document}
\title{The Poset Metrics\\
       That Allow Binary Codes of Codimension $m$\\
       to be $m$-, $(m-1)$-, or $(m-2)$-Perfect%
     \thanks{This is author's version of the correspondence in the
\href{http://ieeexplore.ieee.org/xpl/RecentIssue.jsp?punumber=18}{IEEE Transactions on Information Theory} 54(11) 2008, 5241-5246, Digital
Object Identifier \href{http://dx.doi.org/10.1109/TIT.2008.929972}{10.1109/TIT.2008.929972},
\copyright 2008 IEEE. Personal use of
the material of the correspondence is permitted. However, permission to
reprint/republish the material for advertising or promotional purposes or for
creating new collective works for resale or redistribution to servers or lists,
or to reuse any copyrighted component of this work in other works must be
obtained from the IEEE.
            }%
     \thanks{The results of the paper were presented
             at the IEEE International Symposium on Information Theory ISIT2007, Nice, France
            }%
      }
\author{Hyun~Kwang~Kim
   \thanks{H.~K.~Kim is with the Department of Mathematics,
           Pohang University of Science and Technology,
           Pohang 790-784, South Korea
           (e-mail: \mbox{hkkim@postech.ac.kr})
          }
        and
        Denis~S.~Krotov
   \thanks{D.~S.~Krotov is with the Sobolev Institute of Mathematics,
           prosp. Akademika Koptyuga 4, Novosibirsk, 630090, Russia
           (e-mail: \mbox{krotov@math.nsc.ru})
          }
       }
\maketitle
\begin{abstract}
A  binary poset code of codimension $m$
(of cardinality $2^{n-m}$, where $n$ is the code length) can
correct maximum $m$ errors.
All possible poset metrics that allow codes
of codimension $m$ to be \mbox{$m$-,} \mbox{$(m-1)$-,} or $(m-2)$-perfect are
described. Some general conditions on a poset which guarantee
the nonexistence of perfect poset codes are derived; as examples,
we prove the nonexistence of $r$-perfect poset codes for some $r$
in the case of the crown poset
and in the case of the union of disjoint chains.

\end{abstract}

\section{ Introduction}

We study the problem of existence of perfect codes
in poset metric spaces,
which are a generalization of the Hamming metric space,
see \cite{95BGL}.
There are several papers \cite{03AKKK,04HK,05KO}
on the existence of $1$-, $2$-, or $3$-error-correcting poset codes.
The approach of the present work is opposite;
we start to classify posets that admit the existence of perfect codes
correcting as many as possible errors with respect to the code length and dimension,
i.\,e., when the number of errors is close to the code codimension.

As stated by Lemma \ref{rm} below, the codimension $m$ of an $r$-error-correcting
$(n,2^{n-m})$ code cannot be less than $r$.
And the posets that allow binary poset-codes of codimension $m$ to be $m$-perfect
have a simple characterization (Theorem \ref{th0}).

The main results of this work,
stated by Theorem \ref{th1} and Theorem \ref{l-m-2},
are criteria for the existence of
$(m-1)$- and $(m-2)$-perfect $(n,2^{n-m})$ $P$-codes.
The intermediate results formulated as lemmas may also be
useful for the description of other poset structures admitting
perfect poset codes.

Let $P=([n],\preceq)$ be a poset, where $[n]\triangleq\{1,\ldots,n\}$.
A subset $I$ of $[n]$ is called an \emph{ideal}, or \emph{downset} (an \emph{upset}, or \emph{filter}) iff
for each $a\in I$ the relation $b\preceq a$
(respectively, $b\succeq a$)  means $b\in I$.
For $a_1,...,a_i\in P$ denote by ${<}a_1,...,a_i{>}$ or ${<}\{a_1,...,a_i\}{>}$
the \emph{principal} ideal of $\{a_1,...,a_i\}$,
i.e., the minimal ideal that contains $a_1,...,a_i$;
and by ${>}a_1,...,a_i{<}$ or ${>}\{a_1,...,a_i\}{<}$, the minimal
upset that contains $a_1,...,a_i$.

Denote by ${\mathcal{I}}_P^r\subset 2^{[n]}$
the set of all $r$-ideals (i.\,e., ideals of cardinality $r$) of $P$, where
$r\in\{0,1,...,n\}$.

If $S$ is an arbitrary set (poset), then the set of all subsets of $S$
is denoted by $2^S$. The set $2^{[n]}$ will be also denoted as
$F^n$, and we will not distinguish subsets of ${[n]}$ from
their characteristic vectors; for example,
$2^{[5]}\ni \{2,4,5\}=(01011)\in F^5$.

If $\bar x\in 2^{[n]}$, then the \emph{$P$-weight} $w_P(\bar x)$ of $\bar x$ is the
cardinality of ${<}\bar x{>}$. Now, for two elements $\bar x,\bar y\in F^n$ we
can define the \emph{$P$-distance} $d_P(\bar x,\bar y)\triangleq w_p(\bar x+\bar y)$,
where $+$ means the symmetrical difference in terms of subsets of
$[n]$ and the mod $2$ addition in terms of their characteristic
functions.

For $r\in \{0,...,n\}$ we denote by
${\mathcal{B}}_P^r\triangleq\{\bar x\in F^n \, | \, w_P(\bar x)\leq r\}$ the
 \emph{ball} of radius $r$ with center in the all-zero vector $\bar 0$.
A subset $\mathcal{C}$ of $F^n$ is called an \emph{$r$-error-correcting $P$-code}
(\emph{$r$-perfect $P$-code}) iff each element $\bar x$ of $F^n$ has at most one
(respectively, exactly one) representation in the form
$\bar x=\bar c+\bar b$, where $\bar c\in \mathcal{C}$ and $\bar b\in {\mathcal{B}}_P^r$.
In other words, the  balls of radius $r$ centered
in the codewords of an $r$-error-correcting $P$-code $\mathcal{C}$
are mutually disjoint (the \emph{ball-packing condition}) and, if $\mathcal{C}$ is $r$-perfect,
cover all the space $F^n$. As a consequence,
$$ |\mathcal{C}|\leq |F^n|/|{\mathcal{B}}_P^r| $$
(the \emph{ball-packing bound}), where equality is equivalent to the $r$-perfectness of $\mathcal{C}$.

For the rest of the paper we will use the following notations.
Let $\mathcal{C}\subset F^n$ be a $P$-code and $\bar 0\in \mathcal{C}$;
denote \\
\begin{itemize}
    \item $m\triangleq n-\log_2|\mathcal{C}|$,
    \item $P^r\triangleq \bigcup_{I\in {\mathcal{I}}_P^{r}}I \subseteq [n]$,
    \item $u\triangleq |\bigcap_{I\in {\mathcal{I}}_P^{r}}I|$,
    \item $\widetilde P^r\triangleq  P^r\backslash \bigcap_{I\in {\mathcal{I}}_P^{r}}I$
     (studying $r$-perfect codes, we can call $\widetilde P^r$ the ``essential part'' of $P$;
      indeed, the ball ${\mathcal{B}}_P^r$ is the Cartesian product
      of ${\mathcal{B}}_{\widetilde P^r}^{r-u}$ and $2^{ P^r\backslash \widetilde P^r}$),
    \item $\lambda\triangleq | P^r|-r$,
    \item $\max(R)$ denotes the set of maximal elements of a poset $R$,
    \item $\min(R)$ denotes the set of minimal elements of a poset $R$,
    \item $k \triangleq |\max(\widetilde P^r)|$.
\end{itemize}
Note that $u$, $\lambda$, and $k$ depend on $P$ and $r$ though the notations do not
reflect this dependence explicitly.

\section{ $m$-Error-Correcting Poset Codes}

We start with several auxiliary statements. The first one is easy and well known.

\begin{proposition}\label{II'}\mbox{}\\
{\rm a)} Let $0\leq r \leq r'\leq n$ and $I\in {\mathcal{I}}_P^r$; then there exists
$I'\in{\mathcal{I}}_P^{r'}$ such that $I\subseteq I'$.\\
{\rm b)} Let $0\leq r' \leq r\leq n$ and $I\in {\mathcal{I}}_P^r$; then there exists
$I'\in{\mathcal{I}}_P^{r'}$ such that $I'\subseteq I$.
\end{proposition}
\begin{proof}
If $r'=r$, then $I'=I$ in both cases.

a) In the case $r'=r+1$ let $j$ be a minimal element of $P\backslash I$;
then $I'\triangleq I\cup\{j\}$ satisfies the condition.

b) In the case $r'=r-1$ let $j$ be a maximal element of $I$;
then $I'\triangleq I\backslash\{j\}$ satisfies the condition.

The general cases $r'=r\pm t$ are proved by induction.
\end{proof}

\begin{corollary}\label{I'}
For each $r'$ from $0$ to $n$ the set ${\mathcal{I}}_P^{r'}$ is not empty.
\end{corollary}
\begin{proof}
Assigning $r=0$ and
$I=\emptyset$ in Proposition \ref{II'} we get at least one ideal in
${\mathcal{I}}_P^{r'}$.
\end{proof}

\begin{proposition}\label{B2I}
${\mathcal{B}}_P^r=\bigcup_{I\in{\mathcal{I}}_P^r}2^I$.
\end{proposition}

\begin{proof}
Let $\bar x\in
{\mathcal{B}}_P^r$, i.\,e., $w_P(\bar x)=|{<}\bar x{>}|\leq r$.
By Proposition \ref{II'} there exists an ideal $I\in{\mathcal{I}}_P^r$
such that ${<}\bar x{>}\subseteq I$.
So, we have $\bar x\subseteq {<}\bar x{>}
\subseteq I$, and $\bar x\in 2^I$.

Conversely, if $\bar x\in 2^I$ for some $I\in{\mathcal{I}}_P^r$,
then ${<}\bar x{>}\subseteq I$
and $w_P(\bar x)=|{<}\bar x{>}|\leq |I|=r$.
\end{proof}

Since $|{\mathcal{I}}_P^r|\geq 1$ by Corollary \ref{I'}, we immediately obtain 
\begin{corollary}\label{B2r}
$|{\mathcal{B}}_P^r|\geq 2^r$.
\end{corollary}

The following lemma is straightforward from the ball-packing bound and Corollary \ref{B2r}.
\begin{lemma}\label{rm}
If a $(n,2^{n-m})$ $P$-code $\mathcal{C}\subset F^n$ is $r$-error-correcting,
then $r\leq m$.
\end{lemma}

\begin{theorem}[characterization of $m$-error-correcting  $P$-codes] \label{th0}
An $(n,2^{n-m})$ code $\mathcal{C}$ is an $m$-error-correcting $P$-code
if and only if the following two conditions hold:\\
{\rm a)} ${\mathcal{I}}_P^{m}$ contains exactly one ideal $I$;\\
{\rm b)} there is a function $f:2^{P\backslash I}\to 2^{I}$
such that $\mathcal{C}=\{Y\cup f(Y) | Y\in 2^{P\backslash I}\}$, i.\,e.,
the code $\mathcal{C}$ is \emph{systematical} with
\emph{information symbols} ${P\backslash I}$ and
\emph{check symbols} $I$.

Every $m$-error-correcting $P$-code is an $m$-perfect $P$-code.
\end{theorem}
\begin{proof}
We first  show that a) and b) hold for any $m$-error-correcting $P$-code $\mathcal{C}$.
If ${\mathcal{I}}_P^{m}$ contains more than one ideal, then
$|{\mathcal{B}}_P^r|>2^{m}$, and we have a contradiction with the ball-packing condition.
So, ${\mathcal{I}}_P^{m}$ contains exactly one ideal, say, $I$.

If there is no such a function as in b), then there are
two codewords $\bar c_1,\bar c_2\in \mathcal{C}$
that coincide in ${P\backslash I}$.
Then $\bar c_1+\bar c_2\subseteq I$, and
$d_P(\bar c_1,\bar c_2)=|{<}\bar c_1+\bar c_2{>}|\leq |I|=m$;
therefore $\mathcal{C}$ is not $m$-error-correcting.
So, b) is a necessary condition.

Assume a) and b) hold. We show that $\mathcal{C}$ is an $m$-perfect code.
We need to check that
for each $\bar y\in F^n$ there exists a unique
$\bar c\in \mathcal{C}$ such that $d_P(\bar c,\bar y)\leq m$.
Such $\bar c$ can be defined by
$\bar c=\bar y\cap (P\backslash I)\cup f(\bar y\cap (P\backslash I))$.
It is a code vector by the definition of $f$;
and $d_P(\bar c,\bar y)\leq m$
because $\bar c+\bar y\subseteq I$.
The uniqueness follows from the equalities
$|\mathcal{C}|=2^{n-m}=|F^n|/|{\mathcal{B}}_P^r|$.
\end{proof}

\section{ Useful Statements}

\begin{proposition}\label{r-err}
A $P$-code $\mathcal{C}$ is $r$-error-correcting if and only if for each different
$\bar c_1,\bar c_2\in \mathcal{C}$ and each $I', I''\in {\mathcal{I}}_P^r$ we have
$\bar c_1 + \bar c_2 \not \subseteq I' \cup I''$.
\end{proposition}
\begin{proof}
{\it Only if:} Assume that
there exist $\bar c_1,\bar c_2\in \mathcal{C}$ and $I', I''\in {\mathcal{I}}_P^r$
such that $\bar c_1 + \bar c_2 \subseteq I' \cup I''$. Consider the
vector $\bar v \triangleq \bar c_1+ (\bar c_1 + \bar c_2)\cap I'$.
We have that $d_P(\bar v,\bar c_1)=w_P(\bar v + \bar c_1)=w_P((\bar c_1 + \bar c_2)\cap
I')\leq |I'|=r$. On the other hand,
$d_P(\bar v,\bar c_2)=w_P(\bar v + \bar c_2)
 =w_P(\bar c_1 + \bar c_2+(\bar c_1 + \bar c_2)\cap I')
 =w_P((\bar c_1 + \bar c_2)\backslash I') \leq |I''|=r$
 because
 $(\bar c_1 + \bar c_2)\backslash I'\subseteq I''$ by
 assumption.
 So, 
 $\mathcal{C}$ is not $r$-error-correcting.

{\it If:} Let the $P$-code $\mathcal{C}$ be not $r$-error-correcting. Then there
exist two different codewords $\bar c_1,\bar c_2\in \mathcal{C}$ and a
vector $\bar v\in F^n$
such that $d_P(\bar v,\bar c_1)=|{<}\bar v+\bar c_1{>}|\leq r$
and $d_P(\bar v,\bar c_2)=|{<}\bar v+\bar c_2{>}|\leq r$.
By Proposition \ref{B2I} we have that $\bar v + \bar c_1 \subseteq I'$
and $\bar v + \bar c_2 \subseteq I''$ for some $I',I'' \in {\mathcal{I}}_P^r$. Then
$\bar c_1 + \bar c_2 = (\bar v + \bar c_1)+ (\bar v + \bar c_2) \subseteq I'\cup I''$.
\end{proof}

The statement (Corollary \ref{2-cov}) that we will use for
proving the main result can be derived from each of the
following two lemmas.

\begin{lemma}\label{I1I2}
If there is an $r$-error-correcting
$(n,2^{n-m})$ $P$-code, then $|I'\cup I''|\leq m$ for each $I',I''\in {\mathcal{I}}_P^r$.
\end{lemma}
\begin{proof}
Assume $|I'\cup I''| > m$, i.e., $|P\backslash(I'\cup I'')|<n-m$.
Since $|\mathcal{C}|=2^{n-m}$,
there are two different codewords $\bar c_1,\bar c_2\in \mathcal{C}$
that coincide in $P\backslash(I'\cup I'')$.
This contradicts Proposition \ref{r-err}.
\end{proof}

\begin{lemma}\label{I-I'I''}
Suppose there exists a vector $\bar v\in F^n\backslash {\mathcal{B}}_P^r$ such that for
each $I\in{\mathcal{I}}_P^r$ it is true that $\bar v\cup I\subseteq I'\cup I''$
for some $I',I''\in {\mathcal{I}}_P^r$. Then no $r$-perfect $P$-codes exist.
\end{lemma}
\begin{proof}
Assume the contrary,
i.\,e., there exists an $r$-perfect $P$-code $\mathcal{C}$ and
$\bar 0\in \mathcal{C}$.
Let $\bar c$ be a codeword such that $d_P(\bar v,\bar c)\leq r$.
Then $\bar v + \bar c \in {\mathcal{B}}_P^r$
and by Proposition \ref{B2I} it is true that
$\bar v + \bar c \subseteq I$ for some $I\in{\mathcal{I}}_P^r$.
Therefore $\bar c \subseteq \bar v \cup I$.
By hypothesis,
$\bar c \subseteq I'\cup I''$
for some $I',I''\in {\mathcal{I}}_P^r$, and we get a contradiction with
Proposition \ref{r-err}.
\end{proof}

The following two corollaries are weaker than Lemma \ref{I-I'I''},
but their conditions are more handy for verification.
Given an ideal $V$, denote
 $$W(V)\triangleq [n] \setminus {>}\max(V){<}.$$
It is clear that $W(V)$ is an ideal and it includes $V \setminus \max(V)$.

\begin{corollary}\label{I-I'}
Suppose $V\in{\mathcal{I}}_P^{r+1}$.
Then the following conditions are equivalent and imply
the nonexistence of $r$-perfect $P$-codes:\\
{\rm a)} every $I\in{\mathcal{I}}_P^r$ contains
at least one element $b$ of $\max(V)$;\\
{\rm b)} $|W(V)|<r$.
\end{corollary}
\begin{proof}
Assume a) does not hold, i.\,e.,
there is $I\in{\mathcal{I}}_P^r$ such that $I\cap \max(V)=\emptyset$.
Then $I\subseteq W(V)$ and $|W(V)|\geq |I|\geq r$.
So, b) implies a).

Assume b) does not hold, i.\,e., $|W(V)|\geq r$.
By Proposition \ref{II'} there exists an ideal $I\subseteq W(V)$ of cardinality $r$.
Then $I\cap \max(V)=\emptyset$ and a) does not hold too.
So, a) implies b).

Assume a) holds.
By Lemma \ref{I-I'I''} with $\bar v = V$, $I'=V\backslash \{b\}$, $I''=I$,
we get the nonexistence of $r$-perfect $P$-codes.
\end{proof}

\begin{figure}
\mbox{}\hfill
 \parbox{1,50in}
 {
 \includegraphics[scale=1.0]{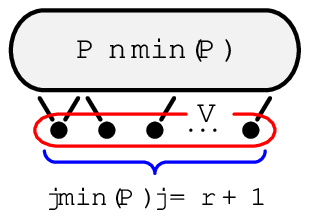}
 \caption{Example \ref{ex1}}\label{fig1}
 }
 \hfill
 \parbox{3,87cm}
 {
 \includegraphics[scale=1.0]{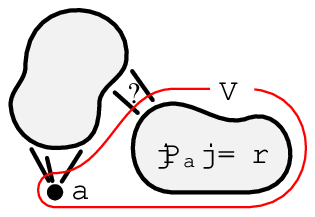}
      \caption{\hfill \mbox{Example \ref{ex2}}}\label{fig2}
 }
\hfill\mbox{}
\end{figure}
\begin{example} \label{ex1}(Fig.~\ref{fig1}.)
If $|\min(P)|=r+1$, then
the ideal $V \triangleq \min(P)$
satisfies the conditions of Corollary \ref{I-I'};
hence, there exist no $r$-perfect $P$-codes.
\end{example}

\begin{example}\label{ex2}(Fig.~\ref{fig2}.)
Let $a$ be a minimal element of $P$ and let $P_a$ be the ideal $P_a \triangleq \{b\,|\,b\not\succeq a\}$.
If $|P_a|=r$, then no $r$-perfect $P$-codes exist,
because
the ideal $V \triangleq \{a\}\cup P_a$
satisfies the conditions of Corollary \ref{I-I'}.
\end{example}

\begin{example}\label{ex3}
Let $P$ be a poset that consists of $t\geq 2$ disjoint chains
and $t-1\leq r < n$.
Then no $r$-perfect $P$-codes exist.
Indeed, it is easy to see that an arbitrary $(r+1)$-ideal $V$
that contains all $t$ minimal elements of $P$
satisfies the conditions of Corollary \ref{I-I'}.
\end{example}

The subcase of Example \ref{ex3}
where $t=2$ and the chains are equipotent
coincides with the binary case of \cite[Theorem 2.2]{95BGL}
(which was proved for codes over arbitrary finite field).

In the next example we see that $r$-perfect $P$-codes do not
exist for sufficiently large $r$ if $P$ is the {\em crown}, i.\,e.,
$n=2t\geq 4$, $i\preceq t+i$, $i+1\preceq t+i$,
$1\preceq 2t$, $t\preceq 2t$,
and these are the only strict comparabilities in $P$.
The existence of $1$-, $2$-, and $3$-perfect crown-codes has been studied in
\cite{03AKKK,05KO}.

\begin{figure}
 \mbox{}\hfill
 \parbox{6,00cm}
 {
 \includegraphics[scale=0.7]{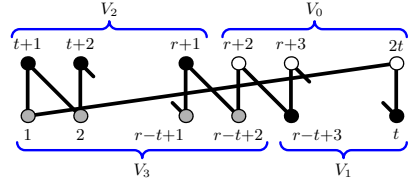}
 \caption{Example \ref{ex4}: the crown poset}\label{fig4}
 }
 \hfill\mbox{}
\end{figure}

\begin{example}\label{ex4}
Let $P$ be a {crown} with $n=2t\geq 6$ and let $t/2\leq r<n$;
then no $r$-perfect $P$-codes exist unless $t=3$ and $r=4$.
Indeed,
it is not difficult to check that condition b) of Corollary~\ref{I-I'}
is satisfied with the following choice of $V$: if $t/2 \leq r < t$
then $V = [t] \backslash \{2,4,\ldots, 2(t-r-1)\}$;
if $t\leq r<n$ then $V=[r+1]$ (Fig.~\ref{fig4}, where $V=V_1+V_2+V_3$).
\end{example}

\begin{corollary}\label{2-cov}
Suppose there exist two different ideals $I',I''\in {\mathcal{I}}_P^r$
such that $ P^r=I'\cup I''$. Then no $r$-perfect $P$-codes exist.
\end{corollary}
\begin{proof}
\emph{Approach 1}: apply Lemma \ref{I-I'I''} with $\bar v= P^r$.
\emph{Approach 2}: if $r$-perfect $P$-codes exist,
then $| P^r|> m$
(indeed, $2^{ P^r}$ includes the ball ${\mathcal{B}}_P^r$ of cardinality $2^m$, but at least one point,
$ P^r$, of $2^{ P^r}$ does not belong to ${\mathcal{B}}_P^r$)
and we have
a contradiction with Lemma \ref{I1I2}.
\end{proof}

Let ${\mathcal{W}}_n^r$ be the set of all $r$-subsets of $[n]$.
Define the distance $d_J(I,I')\triangleq |I+ I'|/2$
(the {\em Johnson distance}) between any $I$ and $I'$ from ${\mathcal{W}}_n^r$.
Let $G_n^r$ be the {\it adjacency graph} of ${\mathcal{W}}_n^r$, where
two subsets $I,I'\in {\mathcal{W}}_n^r$ are {adjacent} iff
$d_J(I,I')=1$; and let $G_n^r(P)$ be the subgraph of $G_n^r$
induced by ${\mathcal{I}}_P^r$.
The following known fact can be easily proved by induction.

\begin{proposition}\label{circ}
Let $I$ and $I'$ be ideals from ${\mathcal{I}}_P^r$.
Then $I$ and $I'$ are connected by a path
of length $d_J(I,I')$ in the graph $G_n^r(P)$.
\end{proposition}
\begin{proof}
If $d_J(I,I')=0\mbox{ or }1$, it is trivial.

Assume the statement holds for $d_J(I,I')=\delta-1\geq 1$.

Let $d_J(I,I')=\delta$.
Let $v$ be a minimal element of $I\backslash I'$ and
$v'$ be a maximal element of $I'\backslash I$.
It is not difficult to check that the
set $I''\triangleq \{v\}\cup I' \backslash \{v'\}$ is an ideal
from ${\mathcal{I}}_P^r$. Since $d_J(I,I'')=\delta-1$ and $d_J(I'',I')=1$,
the induction assumption proves the statement.
\end{proof}

\begin{corollary}\label{conn}
The graph $G_n^r(P)$ is connected.
\end{corollary}

\begin{proposition}\label{ABC}
There exist sequences $I_0,I_1,\ldots,I_{\lambda}\in {\mathcal{I}}_P^r$ and $a_1,\ldots,a_{\lambda}\in [n]$ such
that
\begin{equation}\label{abc}
I_0 \cup I_1 \cup \ldots \cup I_i = I_0 \cup \{a_1, \ldots, a_i\}, \qquad i=0,1,\ldots,\lambda.
\end{equation}
\end{proposition}
\begin{proof}
We construct the sequences by induction.
By Corollary \ref{I'}  there exists $I_0\in {\mathcal{I}}_P^r$.

Assume that for $l<\lambda$ there exist
$I_0,I_1,\ldots,I_{l}\in {\mathcal{I}}_P^r$ and $a_1,\ldots,a_{l}\in [n]$
such that (\ref{abc}) holds for all $i=0,1,\ldots,l$.
We want to find appropriate $I_{l+1}$ and $a_{l+1}$.
Let
$N\triangleq I_0\cup\{a_1,\ldots,a_{l}\}$ and
${\mathcal{I}'}\triangleq \{I\in{\mathcal{I}}_P^r|I\subseteq N\}$.
Since $|N|<| P^r|$, the set ${\mathcal{I}}_P^r\backslash{\mathcal{I}'}$ is not empty.
By Corollary \ref{conn} there are two ideals
$I\in {\mathcal{I}'}$ and $I_{l+1}\in {\mathcal{I}}_P^r\backslash{\mathcal{I}'}$
such that $d_J(I,I_{l+1})=1$.
Then $I_{l+1}$ contains exactly one element from
$ P^r\backslash N$; denote this element by $a_{l+1}$.
So, (\ref{abc}) holds automatically for $i=l+1$.
\end{proof}

We use the last proposition to derive the following bound:
\begin{proposition}\label{B>rl}
$|{\mathcal{B}}_P^r|\geq 2^{r-1}(2+\lambda)$.
\end{proposition}
\begin{proof}
Let $I_0,I_1,\ldots ,I_{\lambda}\in {\mathcal{I}}_P^r$ and $a_1,\ldots ,a_{\lambda}\in [n]$ be
sequences satisfying (\ref{abc}).

Consider the sets ${\mathcal{J}}_0 \triangleq 2^{I_0}$, ${\mathcal{J}}_i \triangleq \{\bar x\subseteq {I_i}|a_i\in \bar x\}$, $i\in [\lambda]$.
We have that
\begin{itemize}
    \item ${\mathcal{J}}_i\subseteq 2^{I_i}$, $i=0,1,\ldots,\lambda$;
    \item the sets ${\mathcal{J}}_0$, ${\mathcal{J}}_1$, \ldots, ${\mathcal{J}}_{\lambda}$ are pairwise disjoint;
    \item $|{\mathcal{J}}_0|=2^r$, $|{\mathcal{J}}_i|=2^{r-1}$, $i=1,\ldots,\lambda$.
\end{itemize}
So, $|{\mathcal{B}}_P^r|\geq |\bigcup_{i=0}^{\lambda}2^{I_i}|\geq
\sum_{i=0}^{\lambda}|{\mathcal{J}}_i|=2^{r-1}(2+\lambda)$.
\end{proof}

Using the ball-packing bound, we deduce the following:
\begin{corollary}\label{lmr}
If $r$-error-correcting $(n,2^{n-m})$ $P$-codes exist,
then $\lambda \leq 2^{m-r+1}-2$.
\end{corollary}

\begin{lemma}\label{lmr-p}
If there exists an $r$-perfect $(n,2^{n-m})$ $P$-code with $r<m$, then
$$m - r < \lambda \leq 2^{m-r+1}-2.$$
\end{lemma}
\begin{proof}
a) $\lambda \leq 2^{m-r+1}-2$ holds by Corollary \ref{lmr}.

b) By Proposition \ref{B2I} we have ${\mathcal{B}}_P^r\subseteq 2^{ P^r}$.
If $r<m$, then $| P^r|> r$ and $ P^r\not\in {\mathcal{B}}_P^r$.
Hence $|{\mathcal{B}}_P^r| < |2^{ P^r}|=2^{r+\lambda}$.
Since $|{\mathcal{B}}_P^r|=2^{m}$ by the ball-packing
condition, we have $2^{m} < 2^{r+\lambda}$, i.\,e., $\lambda > m - r$.
\end{proof}

\section{ $(m-1)$-Perfect Codes}

Applying Lemma \ref{lmr-p} for $r=m-1$ we get the following fact.
\begin{corollary}\label{rm-1} {\rm (Case $r=m-1$.)}
If there exists an $r$-perfect $(n,2^{n-m})$ $P$-code
with $r=m-1$, then $\lambda=2$.
\end{corollary}

In the next proposition we describe the structure of a poset $P$ admitting the existence
of $(m-1)$-perfect $P$-codes. Then, in Proposition~\ref{mu-1exist}, we prove the
existence of $(m-1)$-perfect $P$-codes for admissible posets.
 Theorem~\ref{th1} summarize the results of this section.

\begin{proposition} \label{a1a2a3}
Assume that there exists an $r$-perfect $(n,2^{n-m})$ $P$-code with $r=m-1$.
Then ${\mathcal{I}}_P^r=\{I\cup \{a_1\}, I\cup \{a_2\}, I\cup \{a_3\}\}$,
where $I\in {\mathcal{I}}_P^{r-1}$ and $a_1,a_2,a_3\in [n]\backslash I$.
\end{proposition}
\begin{proof}
By Corollary \ref{rm-1} we have $\lambda=2$. By Proposition
\ref{ABC} there are $I_1,I_2,I_3\in {\mathcal{I}}_P^r$
such that $I_1 \cup I_2 \cup I_3= P^r$.

By Corollary \ref{2-cov} we have that
$I_2 \cup I_3= P^r\backslash \{a_1\}$,
$I_3 \cup I_1= P^r\backslash \{a_2\}$,
and
$I_1 \cup I_2= P^r\backslash \{a_3\}$
for some $a_1,a_2,a_3\in  P^r$.
This implies that
$I_1=I\cup\{a_1\}$,
$I_2=I\cup\{a_2\}$,
and
$I_3=I\cup\{a_3\}$,
where $I=I_1 \cap I_2 \cap I_3$.
It is easy to see that $I$ is an ideal.

By the ball-packing condition we have $|{\mathcal{B}}_P^r|=2^m=2^{r+1}$.
Since $|2^{I_1} \cup 2^{I_2} \cup 2^{I_3}|=2^{r+1}$,
there are no other vectors in ${\mathcal{B}}_P^r$ and there are
no other ideals in ${\mathcal{I}}_P^r$,
i.e., ${\mathcal{I}}_P^r=\{I_1,I_2,I_3\}$.\,\nolinebreak
\end{proof}

\begin{proposition}\label{mu-1exist}
{\rm (Existence of $(m-1)$-perfect $[n,n-m]$ $P$-codes.)}
Let ${\mathcal{I}}_P^{m-1}=\{I\cup \{a_1\}, I\cup \{a_2\}, I\cup \{a_3\}\}$,
where $I\in {\mathcal{I}}_P^{m-2}$ and $a_1,a_2,a_3\in [n]\backslash I$.
Let $\bar h_1,...,\bar h_n\in F^m$.
Assume that $\bar h_i$, $i\in I\cup\{a_1,a_2\}$ are
linearly independent and $\bar h_{a_3}=\sum_{i\in I}\alpha_i\bar h_i+\bar h_{a_1}+\bar h_{a_2}$
where $\alpha_i\in \{0,1\}$, $i\in I$. Then the linear code $\mathcal{C}$
defined by $\mathcal{C}\triangleq\{\bar c\in F^n | \sum_{i\in \bar c} \bar h_i = \bar 0\}$
is an $(m-1)$-perfect $P$-code.
\end{proposition}
\begin{proof}
The $P$-code $\mathcal{C}$ is $(m-1)$-perfect if and only if for each $\bar v\in F^n$ there exists
a unique $\bar e\in {\mathcal{B}}_P^{m-1}$ such that $\bar v + \bar e \in \mathcal{C}$,
i.\,e., $\sum_{i\in \bar v} \bar h_i = \sum_{i\in \bar e} \bar h_i$.
So, it is enough to show that
for each $\bar s\in F^m$ there exists
a unique $\bar e\in {\mathcal{B}}_P^{m-1}$ such that
$\sum_{i\in \bar e} \bar h_i = \bar s$.

Since $\{\bar h_i\}_{i\in I\cup\{a_1,a_2\}}$
is a basis of $F^m$,
for each $\bar s\in F^m$ there exists a (unique) representation
$$\bar s=\sum_{i\in I}\beta_i\bar h_i+\gamma_1
\bar h_{a_1}+\gamma_2 \bar h_{a_2},\quad \beta_i,\gamma_1,\gamma_2\in\{0,1\}.$$
Since $\sum_{i\in I}\alpha_i \bar h_i+\bar h_{a_1}+\bar h_{a_2}+\bar h_{a_3}=\bar 0$,
we can write
$$\bar s=\sum_{i\in I}\beta_i \bar h_i+\gamma_1 \bar h_{a_1}+\gamma_2 \bar h_{a_2}+
\gamma_1\gamma_2
\left(\sum_{i\in I}\alpha_i \bar h_i
{+}\bar h_{a_1}{+}\bar h_{a_2}{+}\bar h_{a_3}\right).$$
By grouping the terms differently we can rewrite it as follows.
$$\bar s=\sum_{i\in I}\beta'_i \bar h_i+\gamma'_1
\bar h_{a_1}+\gamma'_2 \bar h_{a_2}+\gamma'_3 \bar h_{a_3},$$
where $\beta'_i\in\{0,1\}$ and
$(\gamma'_1,\gamma'_2,\gamma'_3)
\in\{(0,0,0),\linebreak[1](0,0,1),\linebreak[1](0,1,0),\linebreak[1](1,0,0)\}$.
This means that $\bar s=\sum_{i\in \bar e} \bar h_i$ for some
$\bar e\in {\mathcal{B}}_P^{m-1}$.
Since $|{\mathcal{B}}_P^{m-1}|=|F^m|$,
such a representation is unique.
\end{proof}

\begin{theorem}\label{th1}
$(m-1)$-perfect $(n,2^{n-m})$ $P$-codes exist
if and only if there are $I\in {\mathcal{I}}_P^{m-2}$ and $a_1,a_2,a_3\in [n] \setminus I$
such that \\
{\rm a)}
if $i,j\in\{1,2,3\}$, $i\neq j$, then ${<}a_i,a_j{>}=  \{a_i,a_j\}\cup I;$\\
{\rm b)}  for each $a\in [n]\backslash  P^r$ there exists $i\in\{1,2,3\}$
such that $\{a_i\}\cup I\subseteq {<}a{>}$.
\end{theorem}
\begin{proof}
By Proposition \ref{a1a2a3} and Proposition \ref{mu-1exist}
$(m-1)$-perfect $(n,2^{n-m})$ $P$-codes exist if and only if
${\mathcal{I}}_P^{m-1}=\{I\cup \{a_1\}, I\cup \{a_2\}, I\cup \{a_3\}\}$
for some $I\in {\mathcal{I}}_P^{m-2}$ and $a_1,a_2,a_3\in [n] \setminus I$.
It is easy to check that this is equivalent to conditions a) and b).
\end{proof}

\section{ More Facts}

Before dealing with the case $r=m-2$, it will be useful to prove some
more facts. We first show that we can restrict ourselves
to consider only the essential part $\widetilde P^r$ of the poset $P$.

\begin{lemma}\label{ppp}
The following statements are equivalent.\\
{\rm a)} There exists an $r$-perfect $P$-code $C$.\\
{\rm b)} There exists an $r$-perfect $ P^r$-code $C'$.\\
{\rm c)} There exists an $(r-u)$-perfect $\widetilde P^r$-code $C''$ (recall $u=| P^r\backslash \widetilde P^r|$).\\
The cardinalities of the codes satisfy $|C''|=|C'|=2^{|P\backslash P^r|}|C|$.
\end{lemma}
%
%
%
\begin{proof}
a)$\Leftrightarrow$b).
By the definition, a perfect code corresponds to a partition of the space
into the balls centered in the code vectors.
In our case, the ball ${\mathcal{B}}_P^r={\mathcal{B}}_{ P^r}^r$ is included
in the subspace $2^{ P^r}$ of the space $2^P$.
Therefore, $2^{ P^r}$ can be partitioned into translations of the balls
if and only if $2^P$ can.

b)$\Leftrightarrow$c).
It is not difficult to see that
${\mathcal{B}}_{\widetilde P^r}^{r-u}={\mathcal{B}}_{ P^r}^r\cap 2^{\widetilde P^r}$
and ${\mathcal{B}}_{ P^r}^r= {\mathcal{B}}_{\widetilde P^r}^{r-u} \times 2^{ P^r \backslash \widetilde P^r}$.
So, if translations of ${\mathcal{B}}_{ P^r}^r$ partition $2^{ P^r}$,
then the intersections with $2^{\widetilde P^r}$ give a partition of $2^{\widetilde P^r}$
into translations of ${\mathcal{B}}_{\widetilde P^r}^{r-u}$.
And vice versa, having a partition of $2^{\widetilde P^r}$ and multiplying it by
$2^{ P^r \backslash \widetilde P^r}$ we get a partition of $2^{ P^r}$.

The relations between the cardinalities immediately follows.
\end{proof}
\begin{lemma}\label{height}
If there is an $r$-error-correcting $(n,2^{n-m})$ $P$-code, then
the height of $\widetilde P^r$ is not more than $m-r$
(the \emph{height} is the maximum length of a chain
in the poset).
\end{lemma}
\begin{proof}
Assume the contrary, i.\,e., $\widetilde P^r$ contains $m-r+1$ pairwise
comparable elements $a_0\preceq a_1\preceq\ldots\preceq a_{m-r}$.
Since $a_0\in \widetilde P^r$,
there exists an ideal $I_1\in{\mathcal{I}}_P^r$ such that $a_0\not\in I_1$.
Then $a_0, a_1,\ldots, a_{m-r}\not\in I_1$.
Since $a_{m-r}\in \widetilde P^r$,
there exists another ideal $I_2\in{\mathcal{I}}_P^r$ such that $a_{m-r}\in I_2$.
Then $a_0, a_1,\ldots, a_{m-r}\in I_2$.
We have that
$|I_1\cup I_2|\geq |I_1\cup\{a_0, a_1,\ldots,a_{m-r}\}|=r+(m-r+1)=m+1>m$,
which contradicts Lemma \ref{I1I2}.
\end{proof}

\begin{proposition}\label{aabc}
Let $U$ be an upset of $P$, $l=|P\backslash U|\leq r$, and
there be an $r$-error-correcting $(n,2^{n-m})$ $P$-code $C$.
Then\\
{\rm a)} $|{\mathcal{B}}_{U}^{r-l}|\leq 2^{m-l}$;\\
{\rm b)} if $|{\mathcal{B}}_{U}^{r-l}|= 2^{m-l}$, then $C$ is $r$-perfect
and $\widetilde P^r\subseteq U$.
\end{proposition}
\begin{proof}
a) It is easy to see that ${\mathcal{B}}_{U}^{r-l}\times 2^{P\backslash U}\subseteq
{\mathcal{B}}_{P}^{r}$.
Since $|{\mathcal{B}}_{P}^{r}|\leq 2^m$ and $|2^{P\backslash U}|=2^l$,
we have $|{\mathcal{B}}_{U}^{r-l}|\leq 2^{m-l}$.

b) If $|{\mathcal{B}}_{U}^{r-l}|= 2^{m-l}$, then $|{\mathcal{B}}_{P}^{r}|= 2^m$ (i.\,e., $C$ is
$r$-perfect)
and ${\mathcal{B}}_{U}^{r-l}\times 2^{P\backslash U}= {\mathcal{B}}_{P}^{r}$.
The last equation means that each $r$-ideal of $P$ is a union of an $(r-l)$-ideal of $U$
and $P\backslash U$, i.\,e., $P\backslash U \subseteq \bigcap_{I\in {\mathcal{I}}_P^r}I$
and, consequently, $\widetilde P^r\subseteq U$.
%
\end{proof}

Recall that $k$ is the number of maximal elements in $\widetilde P^r$ and $\lambda=| P^r|-r$.
(Note that if $|{\mathcal{I}}_P^r|>1$, then $\max( P^r)=\max(\widetilde P^r)$, and thus $k=\max( P^r)$.)

\begin{lemma}\label{l-k}
If there is an $r$-error-correcting $(n,2^{n-m})$ $P$-code and $k\geq \lambda$, then
\begin{equation}\label{k}
  2^{r+\lambda}-2^{r+\lambda-k}\sum_{\sigma=0}^{\lambda-1}\left({k \atop\sigma}\right)\leq 2^m.
\end{equation}
\end{lemma}
\begin{proof}
Every subset of $ P^r$ with not more than $r-(r+\lambda-k)$ elements of
$\max(\widetilde P^r)$ belongs to ${\mathcal{B}}_{P}^{r}$, because its principal ideal
contains the same number of elements of $\max(\widetilde P^r)$ and at most
$| P^r|-|\max(\widetilde P^r)|=(r+\lambda-k)$ other elements (in total, not more than $r$).
So, the number of such subsets does not exceed $|{\mathcal{B}}_{P}^{r}|\leq 2^m$.
On the other hand, this number can be calculated as $|2^{ P^r}|$ minus
the number
$
2^{r+\lambda-k}\sum_{j=k-\lambda+1}^{k}\left({k \atop j}\right)
=
2^{r+\lambda-k}\sum_{\sigma=0}^{\lambda-1}\left({k \atop k-\sigma}\right)
=
2^{r+\lambda-k}\sum_{\sigma=0}^{\lambda-1}\left({k \atop \sigma}\right)
$
of subsets that have more than $k-\lambda$ elements in $\max(\widetilde P^r)$.
\end{proof}

\begin{corollary}\label{c-k}
If $m$, $r$ and $\lambda$ are fixed, then there are only finite
number of values of $k$ that admit the existence
of an $r$-error-correcting $(n,2^{n-m})$ $P$-code.
\end{corollary}
\begin{proof}
It is easy to see that
$2^{r+\lambda}-2^{r+\lambda-k}\sum_{\sigma=0}^{\lambda-1}
\left({k \atop \sigma}\right)
\to 2^{r+\lambda}$
as $k\to\infty$.
By Lemma \ref{lmr-p} if $r<m$, then $r+\lambda>m$ and Lemma \ref{l-k}
proves the statement for $r<m$
(taking into account that only finite number of values of $k$ violate the assumption $k\geq \lambda$).

If $r=m$, then by Theorem \ref{th0} we have $\widetilde P^r=\emptyset$
and $k=0$.
\end{proof}

\begin{proposition}\label{a-b-c}\label{abcc}
Assume that $P=\widetilde P^r$.
Recall that $\lambda= n-r$ in this case.
Then\\
{\rm a)} for each $a\in P$ we have $|P\backslash{{<}{a}{>}}|\geq
\lambda$;\\
{\rm b)} if there exists an $r$-perfect $(n,2^{n-m})$ $P$-code,
then for each $a,a'\in P$ we have $|P\backslash{{<}{a,a'}{>}}|\geq
r+\lambda-m$;\\
{\rm c)} for each $b\in P$ we have $|{{>}{b}{<}}|\leq\lambda$;\\
\end{proposition}
\begin{proof}
a) Since $P=\widetilde P^r$,
an element $a$ belongs to at least one $r$-ideal $I$.
Since ${{<}{a}{>}}\subseteq I$, $|I|=r$, and $|P|=r+\lambda$,
there are at least $\lambda$ elements in $P\backslash{{<}{a}{>}}$.

b) As in p. a), there are $r$-ideals $I\ni a$ and $I' \ni a'$.
Since ${{<}{a,a'}{>}}\subseteq I\cup I'$,
 the statement follows from Lemma \ref{I1I2}.

c) Since $P=\widetilde P^r$, there is at least one $r$-ideal $I$
that does not contain $b$.
Then the upset ${>}{b}{<}$ is disjoint with $I$
and its cardinality does not exceed $|P\backslash I|=\lambda$.
\end{proof}
\section{ The Case $r=m-2$}
\begin{theorem}\label{l-m-2}
An $(m-2)$-perfect $(n,2^{n-m})$ $P$-code exists if and only if
$\widetilde P^{m-2}$ is one of the posets illustrated below:\\
 \mbox{}\hfill
 \includegraphics[scale=0.7]{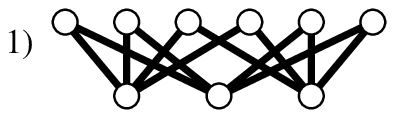} \hfill
 \includegraphics[scale=0.7]{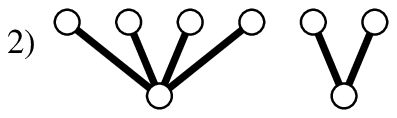} \hfill
 \includegraphics[scale=0.7]{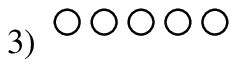}
 \hfill\mbox{}
\end{theorem}
\begin{proof} Assume an $r$-perfect $(n,2^{n-m})$ $P$-code exists with $r=m-2$.
By Lemma \ref{ppp} we can assume $P=\widetilde P^r$.

By Lemma \ref{height} the height of $P$ is $1$ or $2$. So, $P$
consists of maximal and nonmaximal elements, where each nonmaximal
element is also a minimal one. For $a\in P$ denote
$\mbox{valency}(a)\triangleq|\{b\in P\,|\, b\prec a\mbox{ or } b\succ a\}|$.
Proposition \ref{abcc}(a) means that $\mbox{valency}(a)\leq r-1$
for each maximal $a$.  Proposition \ref{abcc}(c) means that
$\mbox{valency}(b)<\lambda$ for each nonmaximal $b$.

By Lemma \ref{lmr-p} we have $\lambda\in\{3,4,5,6\}$.
So, by Corollary \ref{c-k} the number of admissible values
$(\lambda,k)$ is finite. Note that the case $k\leq 2$ is impossible
 by Proposition \ref{abcc}(b).
 The other pairs admitting either (\ref{k}) or $k<\lambda$ are
the following:
$(3,3)$, $(3,4)$, $(3,5)$,
 $(4,3)$, $(4,4)$, $(4,5)$,
  $(5,3)$, $(5,4)$, $(5,5)$, $(5,6)$,
  $(6,3)$, $(6,4)$, $(6,5)$, $(6,6)$.
In all cases we denote by $\{a_1,\ldots,a_k\}$ the set of maximal
elements of $P$. Furthermore, we claim that in all cases except $(3,5)$ the poset contains
at least one nonmaximal element. Indeed, otherwise $|P|=k$, $r=k-\lambda$ and, since $m=r+2$,
we have $|{\mathcal{B}}_P^r|=\sum_{j=0}^{k-\lambda} \left( k \atop j \right) = 2^{k-\lambda+2}$,
which is not true for all considered pairs except $(3,5)$ (in fact, $P=\max(P)$ means
that we have the usual Hamming metric).

{\bf Case $\mathbf{\lambda=3,\ k=3}$.}
%
Let $b_1\in P$ be a nonmaximal element.
W.\,l.\,o.\,g. assume $b_1 \prec a_1$. By Proposition \ref{abcc}(a)
there exists another nonmaximal element $b_2$ which is noncomparable
with $a_1$. W.\,l.\,o.\,g. assume $b_2 \prec a_2$.
By Proposition \ref{abcc}(c)  $\mbox{valency}(b_i)\leq 2$ for $i=1,2$.
The possible cases are:
1) $b_1 {\prec} a_1$, $b_2 {\prec} a_2$;
2) $b_1 {\prec} a_1$, $b_2 {\prec} \{a_2,a_3\}$;
3) $b_1 {\prec} \{a_1,a_2\}$, $b_2 {\prec} a_2$;
4) $b_1 {\prec} \{a_1,a_2\}$, $b_2 {\prec} \{a_2,a_3\}$;
5) $b_1 {\prec} \{a_1,a_3\}$, $b_2 {\prec} a_2$;
6) $b_1 {\prec} \{a_1,a_3\}$, $b_2 {\prec} \{a_2,a_3\}$.
All the cases up to isomorphism  are illustrated in the following figures
(we emphasize that $P$ can have more elements, but in any case it includes an upset
shown in one of the figures, where the dashed lines denotes ``optional'' relations);
Figure (a) corresponds to 1), 3), Figure (b), to 2), 4), 5).
\\[1mm]
 \mbox{}\hfill
 \includegraphics[scale=0.7]{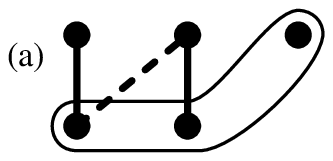}\hfill
 \includegraphics[scale=0.7]{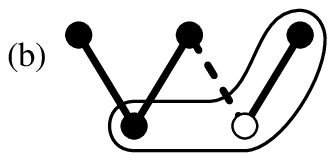}\hfill\mbox{}
\\%
In all these cases we get a contradiction with Corollary~\ref{I-I'} applied for $V$
being the set of all not shown elements of $P$ and those that are banded by
the closed line. The elements of ${>}\max(V){<}$, which
are not in $W(V)$,
are marked by black nodes in the figures.
We see that ${>}\max(V){<}$ has at least $4$ elements; so, $|W(V)|\leq |P|-4 < |P|-\lambda = r$,
and by Corollary~\ref{I-I'}(b) no $r$-perfect $P$-codes exist.

{\bf Case $\mathbf{\lambda=3,\ k=4}$.}
As proved above, $P$ has a nonmaximal element.
Its valency is $1$ or $2$ by Proposition \ref{abcc}(c).
The situation is illustrated by the following figure.
\\[1mm]
 \mbox{}\hfill\includegraphics[scale=0.7]{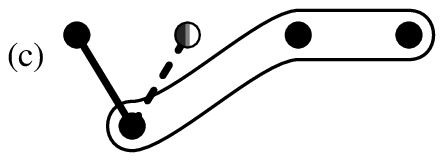}\hfill\mbox{}
\\%
As in the previous case we get a contradiction with Corollary
\ref{I-I'}.

{\bf Case $\mathbf{\lambda=3,\ k=5}$.} By proposition \ref{aabc}(b)
the set of maximal elements coincides with $P$.
In this case
the poset metric coincides with the Hamming metric
and there exists a $2$-perfect \emph{repetition} code
$\{(00000),(11111)\}$.

{\bf Case $\mathbf{\lambda=4,\ k=3}$.} By Proposition \ref{abcc}(b) for
each two maximal elements $a$, $a'$
we have $|P\backslash{{<}{a,a'}{>}}|\geq r+\lambda-m$.
Since $r+\lambda-m=2$, there is a nonmaximal element $b$
noncomparable with $a$ and $a'$. So, there is an upset of $P$
illustrated below
\\[1mm]
 \mbox{}\hfill\includegraphics[scale=0.7]{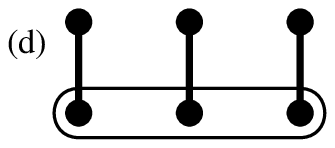}\hfill\mbox{}
\\%
and again we get a contradiction with Corollary \ref{I-I'}.

{\bf Case $\mathbf{\lambda=4,\ k=4}$.}
We claim that there are no subcases different from the ones shown below
\\[1mm]
 \mbox{}
 \includegraphics[scale=0.7]{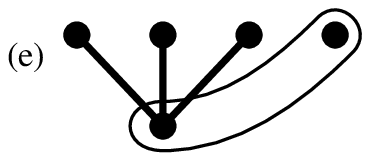}\hfill
 \includegraphics[scale=0.7]{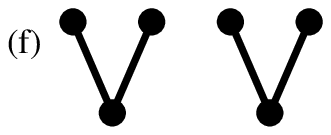}\hfill
 \includegraphics[scale=0.7]{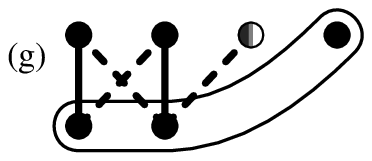}
 \mbox{}
\\%
As in the previous cases there is a nonmaximal element $b_1$ and
$\mbox{valency}(b_1)<4$. Figure (e) illustrates the subcase $\mbox{valency}(b_1)=3$.
Assume $\mbox{valency}(b_1)<3$ and, w.l.o.g., $b_1{\prec}a_1$.
By Proposition \ref{abcc}(a) the set $P\backslash {<}a_1{>}$
contains at least $4$ elements, and one of them,
say $b_2{\prec}a_2$, is not maximal. Figure (f) illustrates
the subcase $b_1{\prec}a_3$, $b_2{\prec}a_4$, or, equivalently,
$b_1{\prec}a_4$, $b_2{\prec}a_3$. Figure (g) illustrates the other
cases.

Subcases (e) and (g) contradict Corollary \ref{I-I'}.
In subcase (f) 
there is at least one more nonmaximal element in $P$,
otherwise $|{\mathcal{B}}_P^r|=12<2^m=16$.
All possibilities to add this
element lead to subcase (e) or (g).

{\bf Case $\mathbf{\lambda=4,\ k=5}$.}
There is nonmaximal element and its valency is 1 (Figure (h)),
2 or 3 (Figure (i)).
\\[1mm]
 \mbox{}\hfill
 \includegraphics[scale=0.7]{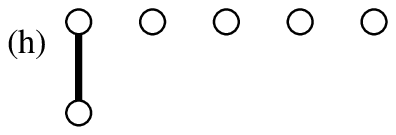}\hfill
 \includegraphics[scale=0.7]{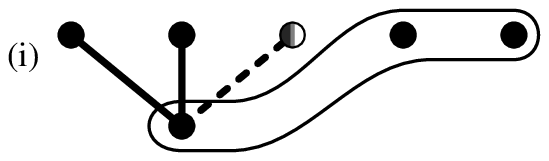}\hfill\mbox{}
\\%
Subcase (h) contradicts Proposition \ref{aabc}(a):
$|{\mathcal{B}}_{U}^{r-l}|=|{\mathcal{B}}_{U}^{2}|=18>16$
(here and below we take $U$ to be the set of all shown elements;
so, $l$ is the number of not shown elements of $P$).
Subcase (i) contradicts Corollary \ref{I-I'}.

{\bf Case $\mathbf{\lambda=5,\ k=3}$.}
Similarly to Case $\lambda=4$, $k=3$, there are
at least two elements of valency $1$ under each $a_i$, $i=1,2,3$ .
We get a contradiction with Corollary \ref{I-I'}, see the figure
\\[1mm]
 \mbox{}\hfill\includegraphics[scale=0.7]{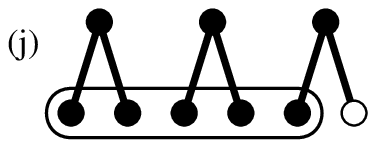}\hfill\mbox{}

{\bf Case $\mathbf{\lambda=5,\ k=4}$.} Let $b_1{\prec} a_1$.
By Proposition \ref{abcc}(a) there is nonmaximal element
noncomparable with $a_1$, say $b_2{\prec} a_2$.
By Proposition \ref{abcc}(b) there is nonmaximal element
noncomparable with $a_1$ and $a_2$.
A contradiction with Corollary \ref{I-I'}, see the figure
\\[1mm]
 \mbox{}\hfill\includegraphics[scale=0.7]{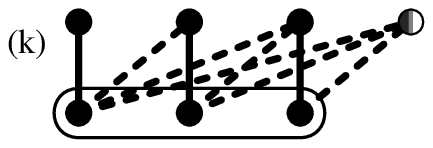}\hfill\mbox{}

{\bf Case $\mathbf{\lambda=5,\ k=5}$.}
Nonmaximal elements cannot be of valency $4$
(Corollary \ref{I-I'}, see figure (l) below).
Assume there is an element $b_1{\prec}a_1$ of valency $1$.
By Proposition \ref{abcc}(a) there is a nonmaximal element
noncomparable with $a_1$, say $b_2$.
Figures (m) and (n) illustrate the cases
$\mbox{valency}(b_2)=1$ and $\mbox{valency}(b_2)\in\{2,3\}$.
\\[1mm]
 \mbox{}\includegraphics[scale=0.7]{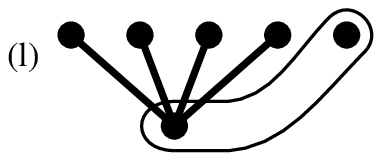}
  \hfill\includegraphics[scale=0.7]{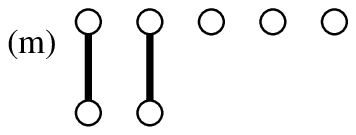}
  \hfill\includegraphics[scale=0.7]{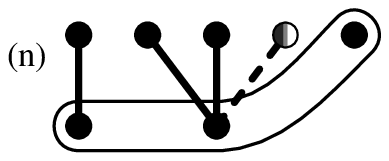}\mbox{}
\\
Case (m) is impossible by Proposition \ref{aabc}(a) with
$|{\mathcal{B}}_{U}^{r-l}|=|{\mathcal{B}}_{U}^{2}|=20>16$.
Case (n) contradicts Corollary \ref{I-I'}.

Let there be a nonmaximal element of valency $2$ or $3$.
By Proposition \ref{abcc}(a) there is another nonmaximal element.
All the situations up to isomorphism are illustrated by the
following figures.
\\[1mm]
 \mbox{}\hfill\includegraphics[scale=0.7]{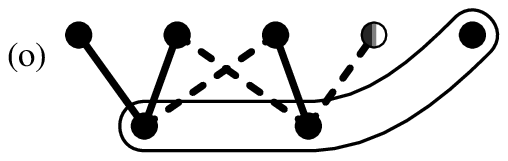}
        \hfill\includegraphics[scale=0.7]{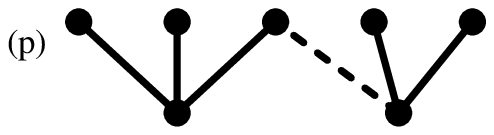}
        \hfill\mbox{}
\\
Case (o) contradicts Corollary \ref{I-I'}.
In the case (p) there is one more element in $P$, otherwise
$|{\mathcal{B}}_{P}|\leq 14<2^m=16$.
But adding a nonmaximal element leads to previous cases.

{\bf Case $\mathbf{\lambda=5,\ k=6}$.} There is a nonmaximal element $b$
with $\mbox{valency}(b)\leq 4$.
The subcase $\mbox{valency}(b)=3$ or $4$ (Figure (q)) is impossible
by Corollary \ref{I-I'};
and the subcase $\mbox{valency}(b)=1$ or $2$ (Figure (r)),
by Proposition \ref{aabc}(a) with
$|{\mathcal{B}}_{U}^{r-l}|=|{\mathcal{B}}_{U}^{2}|>16$.
\\[1mm]
 \mbox{}\hfill\includegraphics[scale=0.7]{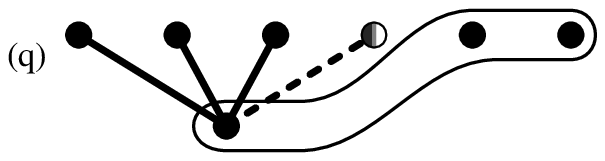}
        \hfill\includegraphics[scale=0.7]{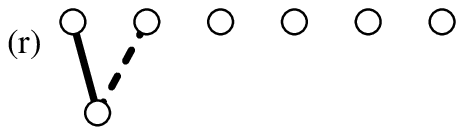}
        \hfill\mbox{}

{\bf Case $\mathbf{\lambda=6,\ k=3}$.}
Similarly to the Cases $(4,3)$ and $(5,4)$, by Proposition \ref{abcc}(b)
every maximal element covers at least three elements of valency $1$.
We have a contradiction with Corollary \ref{I-I'}, see the figure.
\\[1mm]
 \mbox{}\hfill\includegraphics[scale=0.7]{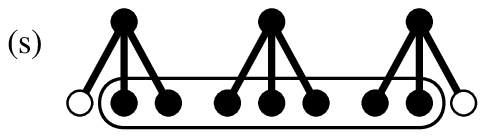}
        \hfill\mbox{}

{\bf Case $\mathbf{\lambda=6,\ k=4}$.}
Let $b_1{\prec}a_1$. By Proposition \ref{abcc}(a) there exists
a nonmaximal $b_2{\not\prec}a_1$, say $b_2{\prec}a_2$.
By Proposition \ref{abcc}(b) there exist
at least two nonmaximal elements noncomparable with $a_1$ and
$a_2$. If there is no more elements in $P$, then $|P|=8$ and
the valency of each maximal element is $1$ by Proposition
\ref{abcc}(a), see Figure (t) below. The other subcase is shown in
Figure (u).
\\[1mm]
 \mbox{}\hfill\includegraphics[scale=0.7]{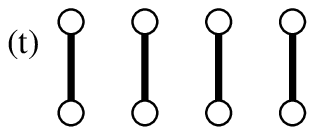}
        \hfill\includegraphics[scale=0.7]{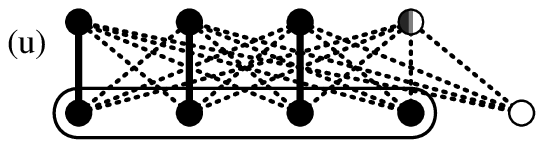}
        \hfill\mbox{}
\\%
In the subcase (t) we have $|{\mathcal{B}}_P|=19$ which is impossible. The
subcase (u) contradicts Corollary \ref{I-I'}.

{\bf Case $\mathbf{\lambda=6,\ k=5}$.}
Let $b_1{\prec}a_1$.
By Proposition \ref{abcc}(a) there exists
a nonmaximal $b_2{\not\prec}a_1$, say $b_2{\prec}a_2$.
By Proposition \ref{abcc}(b) there exists
at least one nonmaximal element noncomparable with $a_1$ and
$a_2$, say $b_3{\prec}a_3$. Let us consider two subcases:
(v) $b_1$, $b_2$, $b_3$ are noncomparable with
$a_4$, $a_5$; (w) otherwise:
\\[1mm]
 \mbox{}\hfill\includegraphics[scale=0.7]{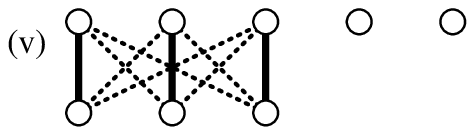}
        \hfill\includegraphics[scale=0.7]{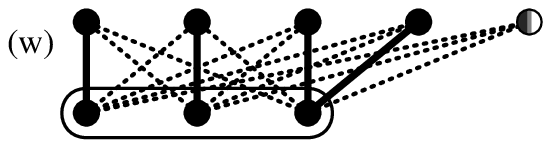}
        \hfill\mbox{}
\\%
Subcase (w) contradicts Corollary \ref{I-I'}.
In Subcase (v), if $\mbox{valency}(a_i)=1$ for some $i\in \{1,2,3\}$,
then we have a contradiction with Proposition \ref{aabc}(a)
($|{\mathcal{B}}_{U}^{r-l}|=|{\mathcal{B}}_{U}^{2}|\geq 18> 2^{m-l}=16$).
Otherwise $|{\mathcal{B}}_{U}^{r-l}|=16$ and by Proposition \ref{aabc}(b)
there are no more elements in $P$. But then there are no $r$-ideals ($r=2$)
that contain $a_1$, $a_2$, or $a_3$, which contradicts our assumption $P=\widetilde P^r$.

{\bf Case $\mathbf{\lambda=6,\ k=6}$.}
The valency of a nonmaximal element is not more than $5$
(Proposition \ref{abcc}(c)). Moreover, the valency $5$ contradicts
  Corollary \ref{I-I'}, see the figure below.
\\[1mm]
 \mbox{}\hfill\includegraphics[scale=0.7]{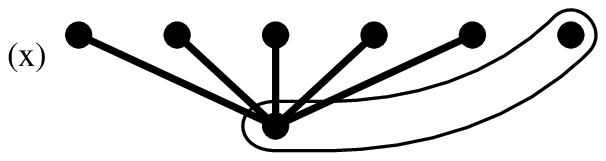}
        \hfill\mbox{}
\\
By Proposition \ref{abcc}(a)
every maximal element is noncomparable with at least one
nonmaximal element.
Calculate the cardinality of the ball ${\mathcal{B}}_P$.
The ideal $P\backslash \{a_1,a_2,a_3,a_4,a_5,a_6\}$ gives
$2^r$ vertices. For each $a_i$ there exists an $r$-ideal containing
only one maximal element $a_i$.
This gives $6\cdot 2^{r-1}$ vertices. So we already have
$2^{r+2}=2^m$ vertices in ${\mathcal{B}}_P$ and, consequently, there is no
other $r$-ideals.
This means that each maximal element is noncomparable with exactly one
nonmaximal element. We already proved that the valency of each
nonmaximal element is not more than $4$, i.\,e., it is noncomparable
with at least two maximal elements. This means that
there are $2$ or $3$ nonmaximal elements. All the cases are
illustrated in the following three figures.
\\[1mm]
  \mbox{}\includegraphics[scale=0.7]{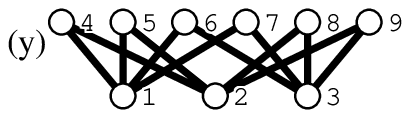}
   \hfill\includegraphics[scale=0.7]{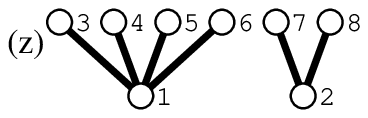}
   \hfill\includegraphics[scale=0.7]{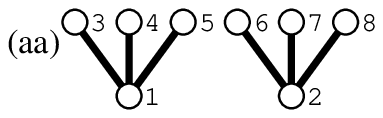}
        \mbox{}
\\
It is known \cite{04HK} that there are perfect
$P$-codes in the cases (y) and (z).
We list here examples of such codes.
$$\begin{array}{rl}
\nonumber 
  (y)&\mbox{linear span}(\{1,2,6,7\},\{1,3,4,5\},\{2,3,8,9\},\{1,4,6,9\}),\\
\nonumber 
  (z)&\mbox{linear span}(\{1,2,3,4\},\{1,2,5,6\},\{1,2,7,8\},\{1,3,5,7\}).
\end{array}$$
The second code is the Hamming $(8,2^4,4)$ code.
The first one is a length $9$ subcode of the Hamming $(16,2^{11},4)$ code.

It remains to show that there are no $2$-perfect $P$-codes in the case (aa).
Assume such a code $C$ exists and, without loss of generality, contains the all-zero vector.
By the definition of perfect $P$-code
$\{2,3\}=\bar c_3+\bar b_3$,
$\{2,4\}=\bar c_4+\bar b_4$, and
$\{2,5\}=\bar c_5+\bar b_5$, where $\bar b_3,\bar b_4,\bar b_5\in {\mathcal{B}}_P^2$ and $\bar c_3,\bar c_4,\bar c_5\in C$.
It is easy to derive from
Proposition~\ref{r-err} that
$\bar b_3\in\{\{4\}, \{1,4\}, \{5\}, \{1,5\}\}$
(indeed, otherwise $\bar c_3$ can be covered by two ideals of size $2$).
So,
$\bar c_3\in \{\{2,3,4\},\{1,2,3,4\},\{2,3,5\},\{1,2,3,5\}\}$; similarly,
$\bar c_4\in \{\{2,3,4\},\{1,2,3,4\},\{2,4,5\},\{1,2,4,5\}\}$,
$\bar c_5\in \{\{2,3,5\},\{1,2,3,5\},\{2,3,5\},\{1,2,3,5\}\}$.
In all cases, two vectors from
$\{2,3,4\}$, $\{1,2,3,4\}$, $\{2,3,5\}$,
$\{1,2,3,5\}$, $\{2,4,5\}$, $\{1,2,4,5\}$ belong to $C$,
and we get a contradiction with Proposition~\ref{r-err}.
\end{proof}

\section*{Acknowledgment}

The authors wish to thank the anonymous referees for their work in
reviewing the manuscript. Their comments have enabled the authors to
greatly improve this correspondence.


\end{document}